\documentstyle[12pt]{article}

\newcommand{\qed}{\nopagebreak\par\noindent\nopagebreak$\blacksquare$\par}
\input amssym.def
\input amssym

\def\intersect{\cap}
\def\Union{\bigcup}
\def\union{\cup}
\def\rmiff{\mbox{ iff }}
\def\rmand{\mbox{ and }}
\def\dom{{\rm dom}}
\def\blocking{{\cal B}}

\def\proof{{\sc Proof. }}
\def\reals{{\Bbb R}}
\def\rationals{{\Bbb Q}}
\def\realfunc{{\reals^\reals}}
\def\almost{{\cal A}}
\def\aa{{{\rm A}(\almost)}}
\def\darboux{{\cal D}}
\def\ad{{{\rm A}(\darboux)}}
\def\A{{\cal H}}
\def\D{{\cal D}}

\def\F{{\cal F}}
\def\G{{\cal G}}
\def\cf{{\rm cf}}
\def\proj{{{\rm pr}_x}}

\def\continuum{{\goth c}}
\def\poset{{\Bbb P}}
\def\posetx{{\Bbb P}^*}
\def\poS{{\Bbb S}}
\def\ONE{{\Bbb I}}
\def\su{\subseteq}
\def\eq{{\goth e}}
\def\la{\langle}
\def\ra{\rangle}
\def\implies{\longrightarrow}

\newtheorem{theorem}{Theorem}[section]
\newtheorem{corollary}[theorem]{Corollary}
\newtheorem{proposition}[theorem]{Proposition}
\newtheorem{lemma}[theorem]{Lemma}
\newtheorem{problem}{Problem}[section]
\newtheorem{example}{Example}[section]
\newtheorem{definition}{Definition}
\newtheorem{remark}{Fact}[section]

\newcommand{\thm}[2]{\begin{theorem}\label{#1}#2\end{theorem}}
\newcommand{\cor}[2]{\begin{corollary}\label{#1}#2\end{corollary}}

\newcommand{\lem}[2]{\begin{lemma}\label{#1}#2\end{lemma}}

\title{Cardinal invariants concerning functions whose sum is
almost continuous.}
\author{}
\date{}

\begin{document}
\maketitle

{\small\noindent Krzysztof Ciesielski\footnotemark[1], Department
of Mathematics, West
Virginia University, Morgantown, WV 26506-6310 (kcies@wvnvms.wvnet.edu)}

\medskip
{\small\noindent Arnold W. Miller\footnote[1]{
The results presented in this paper were
initiated, and partially obtained,
during the Joint US--Polish Workshop in Real Analysis,
{\L}{\'o}d{\'z}, Poland, July 1994.
The Workshop was partially
supported by the NSF grant INT--9401673.
\par
We want to thank Juris Steprans for many helpful conversations.
},
 York University,  Department of Mathematics,
 North York,  Ontario M3J 1P3, Canada, Permanent address:
 University of Wisconsin-Madison,
  Department of Mathematics,
 Van Vleck Hall,
 480 Lincoln Drive,
 Madison, Wisconsin 53706-1388, USA (miller@math.wisc.edu)}

\begin{abstract}
Let $\almost$ stand for the class of all almost continuous functions
from $\reals$ to $\reals$ and let
$\aa$ be the smallest cardinality of a family 
$F\su\realfunc$ for which there is no $g\colon\reals\to\reals$
with the property that $f+g\in\almost$ for all $f\in F$. 
We define cardinal number $\ad$ for the class
$\darboux$ of all real functions with the Darboux property similarly.
It is known, that $\continuum < \aa \leq 2^{\continuum}$ 
\cite{Nat:AC1}.
We will generalize this result by showing that the cofinality of 
$\aa$ is greater that $\continuum$. Moreover, we will show 
that it is pretty much all that can be said about $\aa$
in ZFC, by showing
that $\aa$ can be equal to any regular cardinal between $\continuum^+$
and $2^{\continuum}$ and that it can be equal 
to $2^{\continuum}$ independently of the cofinality of $2^{\continuum}$. 
This solves a problem of T.~Natkaniec \cite[Problem 6.1, p. 495]{Nat:AC1}.

We will also show that $\ad=\aa$ and give a combinatorial characterization 
of this number. This solves another problem of Natkaniec. 
(Private communication.)
\end{abstract}

\section{$\!\!\!\!\!\!\!${\bf.} Preliminaries.}

We will use the following terminology and notation. 
Functions will be identified with their graphs. 
The family of all functions from a set $X$ into $Y$
will be denoted by $Y^X$.
Symbol $|X|$ will stand for the cardinality of a set $X$. 
The cardinality of the set $\reals$ of real numbers is denoted by
$\continuum$.
For a cardinal number $\kappa$ we will write $\cf(\kappa)$
for the cofinality of $\kappa$.
A cardinal number $\kappa$ is regular, if $\kappa=\cf(\kappa)$.
Recall also, that the Continuum Hypothesis (abbreviated as CH)
stands for the statement $\continuum=\aleph_1$.

A function $f\colon\reals\to\reals$ is {\em almost continuous}
(in the sense of Stallings \cite{Stall})
if and only if
for every open set $U\su\reals^2$ containing $f$ 
there exists a continuous function $g\su U$.
So, every neighborhood of $f$ in
the graph topology contains a continuous function.  This
concept was introduced by Stallings \cite{Stall} in connection with
fixed points.
We will use symbol
$\almost$ to denote the family of
almost continuous functions from $\reals$ to $\reals$. 

For $\F\su\reals^\reals$ 
define the cardinal ${\rm A}(\F)$
as follows:
\begin{eqnarray*}
{\rm A}(\F)
& = & \min\{|F|\colon F\su\reals^\reals \&\ \neg\exists
                g\in\reals^\reals\ \forall f\in F\ f+g\in\F\}\\
& = & \min\{|F|\colon F\su\reals^\reals \&\ \forall
                g\in\reals^\reals\ \exists f\in F\ f+g\not\in\F\}
\end{eqnarray*}

For a generalization of the next theorem see Natkaniec \cite{Nat:AC1}.
Fast \cite{fast} proved the same result for the family
of Darboux functions.

\thm{nat}
{$\continuum < \aa \leq 2^{\continuum}$.}\qed

At the Joint US--Polish Workshop in Real Analysis in
{\L}{\'o}d{\'z}, Poland, in July 1994
A.~Maliszewski gave a talk mentioning several problems of his and
T.~Natkaniec. Natkaniec asked whether or not anything more
could be said about the cardinal $\aa$.
(See also Natkaniec \cite[Problem 6.1, p. 495]{Nat:AC1}
or \cite[Problem 1.7.1, p. 55]{Nat:AC2}.)
In what follows we will
show that pretty much nothing more can be said (in ZFC), 
except that the $\cf(\aa)>\continuum$.

We will also study the family
$\darboux\su\realfunc$ of Darboux functions. 
Recall that a function is {\em Darboux} if and only if
it takes every connected set to a connected set, or (in the case
of a real function) satisfies the intermediate value property.
Note that $\almost\su\darboux$. This is because
if for example $f(a)<c<f(b)$ and $c$ is omitted by $f$ on $(a,b)$,
then take the $h$-shape set $H$ (see Figure \ref{fig1}).
\unitlength=1.00mm
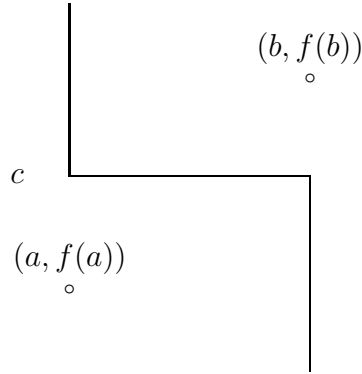
\begin{figure}
\begin{picture}(64.00,64.00)(-32,0)
\put(50.00,7.00){\line(0,1){26.00}}
\put(50.00,33.00){\line(-1,0){32.00}}
\put(18.00,33.00){\line(0,1){23.00}}
\put(18.00,56.00){\line(0,0){0.00}}
\put(18.00,18.00){\circle{1.0}}
\put(50.00,46.00){\circle{1.0}}
\put(50.00,50.00){\makebox(0,0)[cc]{$(b,f(b))$}}
\put(18.00,22.00){\makebox(0,0)[cc]{$(a,f(a))$}}
\put(11.00,33.00){\makebox(0,0)[cc]{$c$}}
\end{picture}
\caption{$h$-shape set $H$\label{fig1}}
\end{figure}
The complement of $H$ is an open neighborhood of the
graph of $f$ which does not contain a graph of a continuous function.
It is known (Stallings \cite{Stall}) that
the inclusion $\almost\su\darboux$
is proper.

It is obvious from the definition that if 
$\F\su\G\su\reals^\reals$ then 
${\rm A}(\F)\leq{\rm A}(\G)$. In particular, 
$\aa\leq\ad$.
At the Joint US--Polish Workshop in Real Analysis in
{\L}{\'o}d{\'z}, Poland, in July 1994,
T.~Natkaniec asked the authors whether it is possible
that ${\rm A}(\almost)<{\rm A}(\darboux)$. We will give a negative answer 
for this question by showing (in ZFC) that 
${\rm A}(\almost)={\rm A}(\darboux)$.

We will finish this section with the following technical fact,
see Natkaniec \cite[Thm. 1.2, p. 464]{Nat:AC1}.

\thm{ThBlock}{(Kellum)
  There exists a family $\blocking$ of closed sets (called {\rm a blocking
  family}) with the properties that:
  \begin{itemize}
    \item for every $f\in\realfunc$ we have
     $$f\in\almost\ \mbox{if and only if }\ \forall B\in\blocking\ 
     f\intersect B\not=\emptyset;$$
    \item for every $B\in\blocking$ the projection $\proj(B)$  of
    $B$ onto the $x$-axis (equivalently, the domain of $B$) 
    is a non-degenerate interval.
  \end{itemize}
 }\qed

The paper is organized as follows. 
We will show that $\ad=\aa$, give some
other characterizations of this cardinal, and prove that
$\cf(\aa)>\continuum$ in Section 2.
In Section 3
we will prove
that some forcing axioms imply that
$\aa$ can be any regular cardinal between $\continuum^+$ and
$2^\continuum$ and that $\aa$ can be equal to $2^\continuum$
for any value of $2^\continuum$. 
The proof of the consistency of the forcing axioms
used in Section 3 will be left for the Section 4.

\section{$\!\!\!\!\!\!\!${\bf.} $\ad=\aa$ and its cofinality.}

We will need the following definitions.

For a cardinal number $\kappa\leq\continuum$
we define the family $$\darboux(\kappa)\su\reals^\reals$$ of
{\em $\kappa$ strongly Darboux functions} as
the family of all functions $f\colon \reals\to\reals$ such that
for all $a,b\in \reals$, $a<b$, and $y\in\reals$ the set
$(a,b)\cap f^{-1}(y)$ has cardinality at least $\kappa$. 

It is obvious from the definition that
\begin{equation}\label{eq1}
\darboux(\lambda)\su
\darboux(\kappa)\ \mbox{ for all cardinals }\ 
\kappa\leq\lambda\leq\continuum.
\end{equation}

We will need the following lemmas. 

\lem{LowerBound}{ ${\rm A}(\darboux(\continuum))>\continuum$.}

\proof Pick a family $F\su \reals^\reals$ of cardinality continuum.
We will find a function $g\in\reals^\reals$ such that
$f+g\in\darboux(\continuum)$ for all $f\in F$. 
Let $$\la\la a_\xi,b_\xi,y_\xi,f_\xi\ra\colon \xi<\continuum\ra$$
be an enumeration of the  set of all
$$\la a,b,y,f\ra\in\reals\times\reals\times\reals\times F\mbox{ with $a<b$},$$
such that each four-tuple appears in the sequence continuum many times.
Define by induction a sequence
$\la x_\xi\in\reals\colon\;\; \xi<\continuum\ra$
such that $$x_\xi\in(a_\xi,b_\xi)\setminus\{x_\zeta\colon\zeta<\xi\}.$$
Then, any function $g\in\reals^\reals$ such that
$g(x_\xi)=y_\xi-f(x_\xi)$ for all $\xi<\continuum$
has the property that
$f+g\in\darboux(\continuum)$ for all $f\in F$. 
\qed

\lem{lemStrDar1}{ $\ad={\rm A}(\darboux(\omega_1))$. }

\proof Since $\darboux(\omega_1)\su\darboux$
we have ${\rm A}(\darboux(\omega_1))\leq\ad$.
To prove the other inequality
let $\kappa={\rm A}(\darboux(\omega_1))$.
Then, by (\ref{eq1}) and Lemma \ref{LowerBound}, 
$$\kappa={\rm A}(\darboux(\omega_1))\geq
{\rm A}(\darboux(\continuum))>\continuum.$$
We will show that $\kappa\geq\ad$.

Let $F\su\reals^\reals$ be a family of cardinality $\kappa$
witnessing $\kappa={\rm A}(\darboux(\omega_1))$:
\begin{equation}\label{eqA}
\forall g\in\reals^\reals\ \exists f\in F\
f+g\not\in\darboux(\omega_1).
\end{equation}
It is enough to find family
$F^*\su\reals^\reals$ of cardinality $\kappa$ such that
\begin{equation}\label{eq2}
\forall g\in\reals^\reals\ \exists f^*\in F^*\
f^*+g\not\in\darboux.
\end{equation}
Define $F^*=\{h\in\reals^\reals\colon \;\exists f\in F\; h=^*f\}$,
where $h=^*f$ if and only if the set $\{x\colon h(x)\neq f(x)\}$
is at most countable.
Since $\kappa>\continuum$ and for every $f\in\reals^\reals$ the set 
$\{h\in\reals^\reals\colon h=^*f\}$ has cardinality $\continuum$, we have
$|F^*|=\kappa$. It is enough to show that $F^*$ satisfies (\ref{eq2}).
So, choose  $g\in\reals^\reals$. Then, by (\ref{eqA}),
there exists $f\in F$ such that
$f+g\not\in\darboux(\omega_1)$.
This means, that there are $a<b$ and $y\in\reals$ such that the set 
$(a,b)\cap (f+g)^{-1}(y)$ is at most countable. 
Then we can find $f^*=^*f$ such that
\begin{itemize}
  \item $(f^*+g)(a)<y$,
  \item $(f^*+g)(b)>y$, and
  \item $(f^*+g)(x)\neq y$ for every $x\in (a,b)$.
\end{itemize}
Thus, $f^*+g\not\in\darboux$. \qed

Now, we are ready for one of our main theorems. 

\thm{equqtion}{ $\ad=\aa$.}

\proof We already know that $\aa\leq\ad$.
So, by Lemma \ref{lemStrDar1},
it is enough to prove that ${\rm A}(\darboux(\omega_1))\leq\aa$.

So, let $\kappa=\aa$. Then, by Theorem \ref{nat},
$\kappa>\continuum$ and,
by the definition of $\aa$, there exists a family $F\su\reals^\reals$
of cardinality $\kappa$
witnessing it, i.e., such that 
\[
\forall g\in\reals^\reals\ \exists f\in F\ f+g\not\in\almost.
\]
In particular, by the definition of the family $\blocking$
of blocking sets (from Theorem \ref{ThBlock}),
\begin{equation}\label{eq3}
\forall g\in\reals^\reals\ \exists f\in F\ 
\exists B\in\blocking\ (f+g)\cap B=\emptyset.
\end{equation}
It is enough to find a family 
$F^*\su\reals^\reals$ of cardinality $\kappa$ such that
\begin{equation}\label{eq4}
\forall g\in\reals^\reals\ \exists f^*\in F^*\ 
f^*+g\not\in\darboux(\omega_1).
\end{equation}
In order to do this, choose a function $h_B\in\reals^\reals$
for every $B\in\blocking$  such that
\[
(x,h_B(x))\in B\ \mbox{ for every }\ x\in\proj(B).
\]
Let 
\[
F^*=\{f-h_B\colon f\in F\ \&\ B\in\blocking\}. 
\]
Clearly
$F^*$ has cardinality $\kappa$, since 
$|\blocking|\leq\continuum<\kappa$. We will show that
$F^*$ satisfies (\ref{eq4}).  Let
$g\in\reals^\reals$. Then, by (\ref{eq3}),
there exist $f\in F$ and $B\in\blocking$ such that
$(f+g)\cap B=\emptyset$. In particular,
\[
[(f-h_B)+g]\cap(B-h_B)=[(f+g)\cap B]-h_B=\emptyset,
\]
where
we define $Z-h_B=\{(x,y-h_B(x))\colon (x,y)\in Z\}$
for any $Z\su\reals^2$.
But $(B-h_B)\supset\proj(B)\times\{0\}$.
Hence,
$[(f-h_B)+g]\cap[\proj(B)\times\{0\}]=\emptyset$. 
In particular, $[(f-h_B)+g]^{-1}(0)\cap\proj(B)=\emptyset$. 
So, $f^*=f-h_B\in F^*$, while $(f-h_B)+g\not\in\darboux(\omega_1)$
since, by Theorem \ref{ThBlock}, $\proj(B)$ contains a non-degenerate 
interval. \qed

To prove the next theorem we need a few more definitions.
For a set $X\su\reals$ and a cardinal number $\kappa\leq\continuum$
we define the family $$\darboux(X,\kappa)\su\reals^X$$
as the family of all functions $f\colon X\to\reals$ such that 
for all $a,b\in X$, $a<b$, and $y\in\reals$ the set
$(a,b)\cap f^{-1}(y)$ has cardinality at least $\kappa$. 
Similarly,
define the cardinal ${\rm A}(\F)$
as before:
\begin{eqnarray*}
{\rm A}(\F)
& = & \min\{|F|\colon F\su\reals^X \&\ \forall
                g\in\reals^X\ \exists f\in F\ f+g\not\in\F\}
\end{eqnarray*}
(Thus $\darboux(\reals,\kappa)=\darboux(\kappa)$.)
It is obvious from the definitions that for $\kappa$ with
$\omega_1\leq\kappa\leq\continuum$
\begin{equation}\label{eqBB}
{\rm A}(\darboux(\reals\setminus\rationals,\kappa))=
{\rm A}(\darboux(\reals,\kappa))
\end{equation}
and also
\begin{equation}\label{eqCC}
{\rm A}(\darboux(X,\kappa))={\rm A}(\darboux(Y,\kappa))\ 
\mbox{ for all order isomorphic $X,Y\su\reals$.}
\end{equation}

\thm{thCof}{ $\aa=\ad={\rm A}(\darboux(\continuum))$.}

\proof By (\ref{eq1}) it is obvious that
$\ad={\rm A}(\darboux(\omega_1))\geq
{\rm A}(\darboux(\continuum))$.

To prove the other inequality let $F\in\reals^\reals$
be a family of cardinality $\kappa$ with $\kappa<\ad.$
It is enough to find $g\in\reals^\reals$ such that
\begin{equation}\label{eqDD}
f+g\in\darboux(\continuum)\ \mbox{ for every }\ f\in F.
\end{equation}
So, let $\la S_\alpha\colon \alpha<\continuum\ra$
be a sequence of pairwise disjoint dense subsets of $\reals$
each of which is order isomorphic to the set
$\reals\setminus\rationals$ of all irrational numbers.
By (\ref{eqBB}) and (\ref{eqCC}) for every $\alpha<\continuum$
we have
\[
\kappa<\ad={\rm A}(\darboux(\omega_1))=
{\rm A}(\darboux(\reals,\omega_1))=
{\rm A}(\darboux(\reals\setminus\rationals,\omega_1))=
{\rm A}(\darboux(S_\alpha,\omega_1)).
\]
We can apply the
definition of ${\rm A}(\darboux(S_\alpha,\omega_1))$
to the family
$$F|_{S_\alpha}=\{f|_{S_\alpha}\in\reals^{S_\alpha}\colon f\in F\}$$
to find a function $g_\alpha\colon S_\alpha\to\reals$
such that
\[
(f|_{S_\alpha})+g_\alpha\in\darboux(S_\alpha,\omega_1)\
\mbox{ for every }\ f\in F.
\]
It is easy to see that any $g\in\reals^\reals$
extending $\bigcup_{\alpha<\continuum}g_\alpha$
satisfies (\ref{eqDD}).

\qed

We will finish this section with one more cardinal
equal to $\aa$. For any infinite cardinal $\kappa$ let
\[
\eq_\kappa=\min\{|F|\colon F\su \kappa^\kappa \&\ \forall
                g\in \kappa^\kappa\ \exists f\in F\ |f\cap g|<\kappa\}.
\]
This cardinal was extensively studied in Landver \cite{land}.

\thm{Comb}{
$\aa=\ad={\rm A}(\darboux(\continuum))=\eq_\continuum$.}

\proof
It is enough to prove that
${\rm A}(\darboux(\continuum))=\eq_\continuum$.
It is also clear that
\[
\eq_\continuum=\min\{|F|\colon F\su \reals^\reals \&\ \forall
                g\in \reals^\reals\ \exists f\in F\ |f\cap g|<\continuum\}.
\]

To prove the inequality
${\rm A}(\darboux(\continuum))\leq\eq_\continuum$
let $F\su\reals^\reals$ have cardinality
$\kappa<{\rm A}(\darboux(\continuum))$.
Then, there exists $g\colon\reals\to\reals$ such that
$g-f\in\darboux(\continuum)$ for every $f\in F$.
In particular,
$|(g-f)^{-1}(0)|=\continuum$, i.e., $f(x)=g(x)$ for continuum
many $x\in\reals$. So, $|f\cap g|=\continuum$ for all $f\in F$, i.e., 
$\kappa<\eq_\continuum$.
This proves
${\rm A}(\darboux(\continuum))\leq\eq_\continuum$.

To prove
$\eq_\continuum\leq{\rm A}(\darboux(\continuum))$
take a family
$F\su\reals^\reals$ of cardinality
$\kappa<\eq_\continuum$. We will show that
$\kappa<{\rm A}(\darboux(\continuum))$.

Choose a sequence
$\la S_{a,b}^y\su(a,b)\colon a,b,y\in\reals, \ a<b\ra$
of pairwise disjoint sets of cardinality continuum.
Applying $\kappa<\eq_\continuum$ to the family
$$F_{a,b}^y=\{(y-f)|_{S_{a,b}^y}\colon f\in F\}$$
we can find $g_{a,b}^y\colon S_{a,b}^y\to\reals$
such that $|(y-f)|_{S_{a,b}^y}\cap g_{a,b}^y|=\continuum$
for every $f\in F$. In particular,
$(y-f)(x)= g_{a,b}^y(x)$, i.e., $(f+g_{a,b}^y)(x)=y$
for continuum many $x\in S_{a,b}^y\su(a,b)$.
Now, if we take any $g\in\reals^\reals$ extending
$\bigcup\{g_{a,b}^y\colon a,b,y\in\reals,\ a<b\}$
then
$(f+g)^{-1}(y)\cap(a,b)$ has cardinality continuum
for every $f\in F$ and $a,b,y\in\reals$, $a<b$.
So, $\kappa<{\rm A}(\darboux(\continuum))$. \qed

\cor{}{$\cf(\aa)>\continuum$.}
\proof
It is obvious that $\cf(\eq_\kappa)>\kappa$ since $\kappa$
can be split into $\kappa$ many sets of size $\kappa$.
\qed

\section{$\!\!\!\!\!\!\!${\bf.} Forcing axioms
and the value of $\aa$.}

In this section we will prove the following two theorems.

\thm{thForMain}{
  Let $\lambda\geq\kappa\geq\omega_2$ be cardinals such that
  $\cf(\lambda)>\omega_1$ and
  $\kappa$ is regular.  Then it is relatively consistent with
  ZFC that the Continuum Hypothesis ($\continuum=\aleph_1$) is true,
  $2^{\continuum}=\lambda$, and $\aa=\kappa$.
}

So for example if $2\leq n\leq 17$, then it is consistent that
$$\continuum=\aleph_1< \aa=\aleph_n\leq\aleph_{17}=2^\continuum.$$

\thm{thForMain2}{
  Let $\lambda$ be a cardinal such that
  $\cf(\lambda)>\omega_1$.
  Then it is relatively consistent with
  ZFC that the Continuum Hypothesis ($\continuum=\aleph_1$) holds and
 $\aa=\lambda=2^{\continuum}$.
}

It follows from Theorem \ref{thForMain2} that $\aa$ can be
a singular cardinal, e.g.
$\aa=\aleph_{\omega_2}$ where $\continuum^+=\omega_2$.
We do not know how to get $\aa$ strictly smaller than $2^\continuum$
and singular.

The technique of proof is a variation on the idea of a Generalized
Martin's Axiom (GMA).
In this section we will formulate the forcing axioms and show that
they imply the results. The proof of the consistency of
these
axioms will be left for Section 4.

For a partially ordered set $(\poset,\leq)$
we say that $G\su\poset$ is a $\poset$-filter if
and only if
\begin{itemize}
  \item for all $p,q\in G$ there exists $r\in G$ with
    $r\leq p$ and $r\leq q$, and
  \item for all $p,q\in\poset$ if $p\in G$ and $q\geq p$, then
  $q\in G$.
\end{itemize}
Define $D\su \poset$ to be dense
if and only if for every $p\in\poset$
there exists $q\in D$ with $q\leq p$.

For any cardinal $\kappa$ and poset $\poset$ define
MA$_\kappa(\poset)$ (Martin's Axiom for $\poset$) to be
the statement that for any family $\D$ of dense
subsets of $\poset$ with $|\D|<\kappa$ there
exists a $\poset$-filter $G$ such that $D\cap G\not=\emptyset$
for every $D\in \D$.

>From now on, let $\poset$ be the following partial order
$$\poset=
\{p\;| \; p:X\to\reals, X\su\reals, \rmand |X|<\continuum\}$$
i.e., the  partial function from $\reals$ to $\reals$ of
cardinality less than $\continuum$.
Define
$p\leq q$ if and only if $q\su p$, i.e., $p$ extends $q$ as a
partial function.

\lem{geq}{MA$_\kappa(\poset)$ implies $\aa\geq\kappa$.}

\proof
We know by Theorem \ref{Comb} that $\aa=\eq_\continuum>\continuum$.
Thus, it is enough to prove that MA$_\kappa(\poset)$ implies
$\eq_\continuum\geq\kappa$ for $\kappa>\continuum$.
Note that for any $\poset$-filter
$G$ since any two conditions in $G$ must have a common extension,
$\Union G$ is a partial function from $\reals$ to $\reals$.
Moreover, it is easy to see that for any $x\in\reals$ the set
$$D_x=\{p\in\poset: x\in\dom(p)\}$$
is dense in $\poset$ and that  
$\Union G\colon\reals\to\reals$ for any $\poset$-filter
$G$ intersecting all sets $D_x$. 

Let $\la S_\alpha:\alpha<\continuum \ra$ be
a partition of $\reals$ into pairwise disjoint sets of size
$\continuum$.
Also for any $f\in{\reals^\reals}$
and $\alpha<\continuum$ the set
$$D_{f,\alpha}=\{p\in\poset: \exists x\in(\dom(p)\cap S_\alpha)\;\;
p(x)=f(x)\}$$
is dense in $\poset$.  Given any $F\subseteq {\reals^\reals}$
with $|F|<\kappa$ let
$$\D=\{D_x:x\in\reals\}\cup\{D_{f,\alpha}:
f\in F,\alpha<\continuum\}.$$
Notice that $|\D|=\continuum<\kappa$. 
Applying MA$_\kappa(\poset)$ we can find a $\poset$-filter $G$ such
that $G$ meets every $D\in\D$.  Letting 
$g=\bigcup G\colon\reals\to\reals$ we
see that $|f\cap g|=\continuum$ for every $f\in\F$.
\qed

The proof of Lemma \ref{geq} is a kind of forcing extension of the
inductive argument used in the proof of Theorem \ref{nat}.

Notice also, that Theorem \ref{thForMain2} follows immediately from
Lemma \ref{geq}, Theorem \ref{nat} and the following theorem.

\thm{thForMain2Proof}{
  Let $\lambda$ be a cardinal such that
  $\cf(\lambda)>\omega_1$.
  Then it is relatively consistent with
  ZFC+CH that $2^{\continuum}=\lambda$ and that
  MA$_\lambda(\poset)$ holds.
}

Thus, we have proved Theorem
\ref{thForMain2} modulo Theorem \ref{thForMain2Proof}.
Theorem \ref{thForMain2Proof} will be proved in
Section 4.

Lemma \ref{geq} shows also one inequality of
Theorem \ref{thForMain}.
To prove the reverse inequality we will use a different partial
order $(\posetx,\leq)$.  It is similar to $\poset$ but in
addition has some side conditions.
\[
\posetx=\{(p,E):p\in\poset \rmand E\su\realfunc
                          \mbox{ with } |E|<\continuum\}.
\]
Define the ordering on $\posetx$ by
\begin{eqnarray*}
(p,E)\leq (q,F) & \rmiff & p\leq q\ \rmand \ E\supseteq F\\
& \rmand & \forall x\in\dom(p)\setminus\dom(q)\;\;
\forall f\in F\; p(x)\not=f(x).
\end{eqnarray*}
The idea of the last condition is that we
wish to create a generic function $g\in \realfunc$ with the property
that for many $f$ we have $g(x)\not=f(x)$
for almost all $x$.
Thus, the condition $(q,F)$ `promises' that for
all new $x$ and old $f\in F$ it should be that $g(x)\not=f(x)$.

For a cardinal number $\kappa$ define Lus$_\kappa(\posetx)$
to be the statement:
\begin{quote}
There exists a sequence $\la G_\alpha:\alpha<\kappa\ra$
of $\posetx$-filters,
called a $\kappa$-Lusin sequence,
such that for every dense set $D\su\posetx$
\[
|\{\alpha<\kappa\colon \;G_\alpha\cap D=\emptyset\}|<\kappa.
\]
\end{quote}
Thus we have a Lusin sequence of $\posetx$-filters.
This is also known as a kind of Anti-Martin's Axiom.
See vanDouwen and Fleissner \cite{df},  Miller and Prikry  \cite{mp},
Todorcevic \cite{tod}, and Miller \cite{sur} for a similar axiom.

\lem{leq}{
  Suppose $\continuum<\kappa$, $\kappa$ is regular, and
  Lus$_\kappa(\posetx)$.
  Then $\aa\leq\kappa$.
}

\proof
Let $\la G_\alpha:\alpha<\kappa\ra$ be a $\kappa$-Lusin sequence
of $\posetx$-filters and let
\[
g_\alpha=\Union \{p:\exists F\; (p,F)\in G_\alpha\}.
\]
Then $g_\alpha$ is a partial function from $\reals$ into $\reals$.
Similarly to the last proof, let
$$D_x=\{(p,F)\in\posetx\colon \; x\in\dom(p)\}.$$
To see that $D_x$ is dense let $(q,F)$ be an arbitrary element
of $\posetx$ and suppose it is not already an element of $D_x$.
The set $Q=\{f(x):f\in F\}$ has cardinality less than $\continuum$
so there exists
$y\in\reals\setminus Q$.  Let $p=q\union\{(x,y)\}$. 
Then $(p,F)\leq (q,F)$ and
$(p,F)\in D_x$. Thus, each $D_x$ is dense in $\posetx$. 
Hence, since $\continuum<\kappa$ and $\kappa$ is regular,
we may assume the each $g_\alpha$ is a total function.

For each $f\in\realfunc$ define
$$D(f)=\{(p,E)\in\posetx: \; f\in E\}.$$

Note that for any $(p,F)$ if we let $E=F\union\{f\}$, then
$(p,F)\leq (p,E)$. Hence $D(f)$ is dense.

Next, note that by the nature of definition of $\leq$ in $\posetx$,
if $(p,F)\in G$, where $G$ is a $\posetx$-filter, and
$g=\Union \{p\colon\exists F\;(p,F)\in G\}$,
then for any $f\in F$ we have $g(x)\not=f(x)$
except possibly for the $x$ in
the domain of $p$.  Therefore for any $f\in \realfunc$ there
exists $\alpha<\kappa$ such that
$|g_\alpha\cap f|<\continuum$.
Thus, the family $\{g_\alpha\colon\alpha<\kappa\}$ shows
that $\aa=\eq_\continuum\leq\kappa$ as was
to be shown.
\qed
\lem{equiv}{ For any regular $\kappa$ we
have Lus$_\kappa(\posetx)\implies
$MA$_\kappa(\posetx)\implies$MA$_\kappa(\poset)$.
}
\proof
This first implication needs that $\kappa$ is regular but is
true for any partial order.
Given a family $\D$ of
dense subsets of $\posetx$ of cardinality less than $\kappa$
and $\la G_\alpha:\alpha<\kappa\ra$ a Lusin sequence for
$\posetx$ it must be that for some $\alpha<\kappa$ that
$G_\alpha$ meets every element of $\D$.

The second implication follows from the fact that in some
sense $\poset$ is `living inside' of $\posetx$.
Let $r:\reals\to\reals$ be a map with
of $|r^{-1}(y)|=\continuum$ for every $y\in\reals$.
Define
$$\pi:\posetx\to\poset
\mbox{ by } \pi(p,F)=r\circ p.$$
Notice that if $(p,E)\leq(q,F)$ then  $\pi(p,E)\leq\pi(q,F)$.
This implies that $\pi(G)$ is a $\poset$-filter for any 
$\posetx$-filter $G$.  Furthermore, we claim that
if $D\su\poset$ is dense, then $\pi^{-1}(D)$ is dense in
$\posetx$.  To see this, let $(p,F)\in\posetx$ be arbitrary.
Since $D$ is dense, there exists $q\leq \pi(p,F)$ with $q\in D$.
Now, find $s\in\poset$ extending $p$ such that 
$r\circ s= q\supseteq r\circ p$ and $s(x)\neq f(x)$
for every $x\in\dom(s)\setminus\dom(p)$ and $f\in F$. 
This can be done by choosing 
\[
s(x)\in r^{-1}(q(x))\setminus\{f(x)\colon f\in F\}
\]
for every $x\in\dom(q)\setminus\dom(p)$.
Then, $(s,F)\leq (p,F)$ and $(s,F)\in\pi^{-1}(q)\su\pi^{-1}(D)$.

This gives
us the second implication, since if $\D$ is a family
of dense subsets of $\poset$ with $|\D|<\kappa$
and $G$ is a $\posetx$-filter
meeting each element of $\{\pi^{-1}(D):D\in\D\}$, then
$\pi(G)$ is a $\poset$-filter meeting each element of $\D$.
\qed

It follows from
Lemmas \ref{geq}, \ref{leq}, and \ref{equiv} 
that Lus$_\kappa(\posetx)$ implies $\aa=\kappa$.
In particular, Theorem \ref{thForMain} follows from
the following theorem.

\thm{thForMainProof}{
  Let $\lambda\geq\kappa\geq\omega_2$ be cardinals such that
  $\cf(\lambda)>\omega_1$ and $\kappa$ is regular.
  Then it is relatively consistent with
  ZFC+CH that $2^{\continuum}=\lambda$ and
  Lus$_\kappa(\posetx)$ holds.
}

Theorem \ref{thForMainProof} will be proved in Section 4.

\section{$\!\!\!\!\!\!\!${\bf.}
Consistency of our forcing axioms.}

In this section we will prove Theorems \ref{thForMain2Proof}
and \ref{thForMainProof}.  For Theorem \ref{thForMain2Proof},
start with a model of GCH and extend it by forcing
with the countable partial functions from $\lambda$ to
$\omega_1$.  For Theorem \ref{thForMainProof} start with a
model of $$2^\omega=\omega_1+2^{\omega_1}=\lambda$$ and
do a countable support iteration of $\posetx$ of length
$\kappa$. $\posetx$ is isomorphic to the eventual dominating
partial order.
For the expert this should suffice.  The rest
of this section is included for
our readers who are not set theorists.
For similar proofs see for example Kamo \cite{kam} and
Uchida \cite{uc}.

We begin with some basic forcing terminology and facts.
(See Kunen \cite{Kun}.)
For a model $M$ of set theory ZFC and a partial
order set $(\poS,\leq)$ a filter $G\su\poS$ is
{\em $\poS$-generic over $M$} if $G$ intersects every dense
$D\su\poS$ belonging to $M$. The fundamental theorem of forcing
states that for every model $M$ of ZFC and every partial order
 $\poS$ from $M$ there exists model $M[G]$ of ZFC
(called an {\em $\poS$-generic extension of $M$})
such that $G$ is $\poS$-generic over $M$ and
$M[G]$ is the smallest model of ZFC such that
$M\su M[G]$ and $G\in M[G]$.
Thus, the simplistic idea for getting
MA$_\kappa(\poset)$ is to start with
model $M$ of $ZFC$, take $\poset$ from $M$ and
look at the model $M[G]$, where $G$ is
$\poset$-generic over $M$. Then, $G$ intersects
``all'' dense subsets of $\poset$ and we are done. 
There are, however, two problems with this
simple approach. First, ``all''
dense subsets of $\poset$ means ``all dense subsets from $M$''
and we like to be able to talk about all dense subsets from
our universe, i.e., from $M[G]$. 
Second, our partial order is 
a set described by some formula as the set having some properties.
There is no reason, in general, that the 
same description will give us the same objects in
$M$ and in $M[G]$.

The second problem will not give us much trouble.
For the generic extensions we will consider, the 
definition of $\poset$ will give us the same objects
in all models we will consider. In the case of the
partial order $\posetx$
this will not be the case, but the new
orders $\posetx$ will be close enough to the old
so that it will not bother us.

To take care of the first of the mentioned problems, we will 
be constructing a Lusin sequence
$\la G_\alpha\colon\alpha<\kappa\ra$ by some kind of induction 
on $\alpha<\kappa$: our final model
can be imagined as $N=M[G_0][G_1]\ldots[G_\alpha]\ldots$
and we will make sure that
every dense subset $D\in N$ of $\posetx$ is taken care of
from some stage $\alpha<\kappa$.

We need some more definitions and facts. 
Given a partial order we say that
$p,q$ are {\em compatible} if there exists $r$
such that $r\leq p$ and $r\leq q$.
A partial order is {\em well-met} provided
for any two elements $p,q$ if
$p$ and $q$ are compatible, then they have a greatest lower
bound, i.e., there exists $r$ such that $r\leq p$ and $r\leq q$ and
for any $s$ if $s\leq p$ and $s\leq q$, then $s\leq r$.
Notice that both partial orders $\poset$ and $\posetx$
used in Lemmas \ref{geq} and
\ref{leq} are well-met. For the case of $\posetx$
if $(p,E)$ and $(q,F)$ are compatible,
then $(p\union q,E\union F)$ is there greatest lower bound.
A subset $L$ of a partial order
is {\em linked} if any two elements of
$L$ are compatible.
A partial order is {\em $\omega_1$-linked} provided it is
a union of $\omega_1$ linked subsets.  Assuming the
Continuum Hypothesis note that the poset $\poset$ used in
the proof of Lemma \ref{geq} has cardinality $\omega_1$
hence it is $\omega_1$-linked.  Note that for any $p\in\poset$
if we define
$$L_p=\{(q,F)\in\posetx: q=p\},$$
then $L_p$ is a linked subset of $\posetx$, hence $\posetx$ is
also $\omega_1$-linked.
A subset $A$ of a partial order
is an {\em antichain} if any two elements of
$A$ are incompatible.
We say that a partial order has the 
{\em $\omega_2$-chain condition} ({\em $\omega_2$-cc})
if
every its antichain has
cardinality less than $\omega_2$.  Clearly $\omega_1$-linked implies
the $\omega_2$-chain condition.  Finally we say a  
partial order is
countably closed if
any descending $\omega$-sequence $\la p_n:n\in\omega\ra$
(i.e., $p_{n+1}\leq p_n$ all $n$) has a lower bound. Notice that
both of our partial orders are countably closed.

All partial orders we are going to consider here 
will be countably closed and will satisfy
$\omega_2$-chain condition. In particular, it is known that if
the generic extension $M[G]$ of $M$ is obtained with such 
partial order,
then $M[G]$ and $M$ have the same cardinal numbers,
the same real numbers, the same countable subsets of real numbers and
the same sets $\reals^X$ for any countable set $X\in M$.
In particular, $\poset$ will be the same in 
$M[G]$ as in $M$.

Let us also notice that every dense set contains
a maximal antichain and if $A$ is a maximal 
antichain, then $D=\{p:\exists q\in A\;\; p\leq q\}$
is a dense set.  Thus a filter $G$ is $\poS$-generic
over a model $M$ if and only if 
it meets every maximal antichain in $M$.

\bigskip

{\sc Proof of Theorem \ref{thForMain2Proof}.}
Take a model $M$ of ZFC+GCH. For a set $X$ in $M$ let
\[
\poS_X=
\{p\in\poset^X\colon p(x)=\emptyset
\mbox{ for all but countably many $x\in X$}\}.
\]
Define an ordering on $\poS_X$ by $p\leq q$ if and only if
$p(x)\leq q(x)$ for every $x\in X$.

Now, let $\lambda$ be as in Theorem \ref{thForMain2Proof}
and let $G$ be a $\poS_\lambda$ generic over $M$.
We will show that MA$_\lambda(\poset)$
holds in $M[G]$.

It is easy to see that $\poS_\lambda$ is countably closed.
It is also known that $\poS_\lambda$ satisfies $\omega_2$-cc
and that
$2^{\omega_1}=\lambda$ in $M[G]$. (See
Kunen \cite[Ch. VII, Lemma 6.10 and Thm. 6.17]{Kun}.)

Now, for $\alpha<\lambda$
let $G_\alpha=\{p(\alpha)\colon p\in G\}$. Then, each $G_\alpha$
is a filter in $\poset$. We will show that for every family
$\D$ of dense subsets of $\poset$ with $|\D|<\lambda$ 
there exists $\alpha<\kappa$ such that $G_\alpha$ intersects every $D$ 
from $\D$. 

In order to argue for it we need two more facts about 
forcing $\poS_X$. (See Kunen \cite[Ch. VII]{Kun}: 
Thm. 1.4 and 2.1 for (A) and Lemma 5.6 for (B).)
\begin{description}
\item{(A)} If $X,Y\in M$  are disjoint and $G$ is $\poS_{X\cup Y}$-generic
   over $M$, then $G_X=G\cap\poS_X$ is $\poS_X$-generic
   over $M$, $G_Y$ is $\poS_Y$-generic over $M[G_X]$, and
   $$M[G_X][G_Y]=M[G].$$
\item{(B)} If $A\su M$ then there exists $X\in M$ with
    $|X|\leq |A|+\omega_1$
    such that $A\in M[G_X]$.
\end{description}

Now, let
$G_\lambda$ be $\poS_\lambda$ generic over
$M$ and let $\D\in M[G_\lambda]$ be a family of dense subsets of $\poset$
with $|\D|<\lambda$.  Let $\A$ be a family of maximal antichains,
one contained in each element of $\D$.  
Then, $|A|\leq\omega_1$
for each $A\in \A$, since $\poset$ satisfies $\omega_2$-cc.
So, by (B), there is
$X\su\lambda$ from $M$ of cardinality $|\A| \cdot \omega_1<\lambda$ 
such that $\A\in M[G_X]$.
Choose $\alpha\in\lambda\setminus X$.
Then since $G_\alpha$ is $\poset$-generic over
$M[G_X]$ it follows that
$G$ meets each element of $\A$ hence of $\D$.
\qed

\bigskip

Next we prepare to prove  Theorem \ref{thForMainProof}.
As mentioned in the beginning of the section, we will try to prove it by
defining some sequence
$\la \poS_\alpha\colon\alpha\leq\kappa\ra$
of partial orders and try to obtain our
final model as $N_\kappa=M[G_\kappa]$
where every $G_\alpha$ is an $\poS_\alpha$-generic
over an appropriate initial model. This technique is called
iterated forcing and needs a few words of introduction.

We can define in $M$ an iterated forcing $\la \poS_\alpha:\alpha<\kappa\ra$
by induction on $\alpha$.
At successor stages we define
$$\poS_{\alpha+1}=\poS_\alpha\times{\posetx}^{M[G_\alpha]}.$$
where ${\posetx}^{M[G_\alpha]}$ is $\posetx$ in the sense of
$M[G_\alpha]$.  (Since we add new elements of $\reals^\reals$
the partial order $\posetx$ changes as our models increase.)
We can't really do it
precisely this way, because $\poset_{\alpha+1}$ must
be in $M$.  However, it is possible
to find its approximation, $\hat{\poset}_\alpha$, in $M$, called
a name for $\poset_\alpha$, and use this instead.
(See Kunen \cite[Ch. VII sec. 5]{Kun}).

For limit ordinals $\lambda<\kappa$, define $\poS_\lambda$
to a set of functions $f$ with domain $\lambda$ such
that $f|_{\alpha}\in \poS_\alpha$ for each $\alpha<\lambda$
and $f(\alpha)=\ONE$ for all but countable many $\alpha$.
Here we use $\ONE$ to denote the largest element of any partial
order.  Countable support iterations originated with
Laver \cite{lav}.  For
details see Baumgartner \cite{Baum} or
Kunen \cite[Ch. VII sec. 7]{Kun}.

The proof that follows will involve a basic lemma used
to show various generalizations of Martin's Axiom hold
for one cardinal up. (See Baumgartner \cite{Baum}
and Shelah \cite{sh}). In particular, we will need the following
theorem.

\thm{thBaum}{ (Baumgartner)
 Assume the Continuum Hypothesis.  Suppose
 $\la \poS_\alpha:\;\alpha<\kappa\ra$
 is a countable support iteration of countably closed
 well-met $\omega_1$-linked partial orders.
 Then for every
 $\alpha\leq\kappa$ we have that $\poS_\alpha$ is
 countably closed and satisfies the $\omega_2$-chain condition.
}

Actually we need only a very weak version of this theorem,
for example, something analogous
to \cite[Theorem VII, 7.3]{Kun} of Kunen.

Now, we are ready for the proof of Theorem \ref{thForMainProof}.

{\sc Proof of Theorem \ref{thForMainProof}.}
Take a model $M$ of ZFC+CH in which
$2^{\continuum}=\lambda$, and $\kappa$ is a regular
cardinal with $\omega_2\leq\kappa\leq\lambda$.
Let $\poS_\alpha$
be a countable support iteration
$\{\poset_\alpha\colon\alpha<\kappa\}$, where
$\poset_\alpha={\posetx}^{M[G^\alpha]}$ for all $\alpha<\kappa$.
Here for $\alpha<\kappa$ let $G^\alpha=G^\kappa|_\alpha$.
Then $G^\alpha$ is $\poS_\alpha$-generic
filter over $M$.

Let $G^\kappa$ be an $\poS_\kappa$-generic
filter over $M$. We will show that
Lus$_\kappa(\posetx)$ holds in $M[G^\kappa]$.

In the model $M[G^\alpha]$
the partial order ${\posetx}^{M[G^\alpha]}$ can be decoded
from $\poS_\alpha$
and we can also decode a filter $G_\alpha$ which is 
${\posetx}^{M[G^\alpha]}$-generic over $M[G^\alpha]$.
We claim that the sequence
$\la G_\alpha\colon\alpha\in \kappa\ra$ is a
Lusin sequence for $\posetx$ in $M[G^\kappa]$.

So, let $D\in M[G^\kappa]$ be a dense subset of $\posetx$
and let $A\in M[G]$ be a maximal antichain contained in 
$D\su\posetx$. Then, $|A|\leq\omega_1$, since
$\posetx$ satisfies $\omega_2$-cc.
So, by the fact similar to (B) above, there is
$\beta<\kappa$
such that $A\in M[G^\beta]$. Then, for every 
$\alpha\geq\beta$, the filter $G_\alpha$
is generic over $M[G^\alpha]\supseteq M[G^\beta]$ 
and so, $G_\alpha$ intersects
both $A$ and $D$.
Therefore, the set
\[
\{\alpha<\kappa\colon G_\alpha\cap D=\emptyset\}
=
\{\alpha<\kappa\colon G_\alpha\cap A=\emptyset\}\su\beta
\]
has cardinality less then $\kappa$. \qed
\bigskip

It is worth mentioning that some
generalizations of the these theorems
are possible where the Continuum Hypothesis
fails.

\end{document}